\documentclass[12pt]{article}
\usepackage{}
\usepackage[numbers,sort&compress]{natbib}

\usepackage{color}

\usepackage[english]{babel}
\usepackage{t1enc}
\usepackage{mathrsfs}
\usepackage[english]{babel}
\usepackage{amsmath,amsthm}
\usepackage{amsfonts}
\usepackage{latexsym}
\usepackage{graphicx}
\usepackage{txfonts}
\usepackage{lineno}
\usepackage{tikz}
\usepackage{tikz-cd}
\usepackage{tikz-qtree}
\usepackage{diagbox}

\usepackage{booktabs}
\newcommand{\rank}{\text{rank}}

\newtheorem{theorem}{Theorem}
\newtheorem{corollary}[theorem]{Corollary}

\newtheorem{lemma}{Lemma}
\newtheorem{remark}{Remark}

\newtheorem{example}{Example}
\newtheorem{definition}{Definition}

\begin{document}
\title{A new criterion for oriented graphs to be determined by their generalized skew spectrum	}
	\author{\small Yiquan Chao$^{{\rm a}}$\quad\quad Wei Wang$^{\rm a}$\quad\quad Hao Zhang$^{\rm b}$\thanks{Corresponding author: zhanghaomath@hnu.edu.cn}
		\\
		{\footnotesize$^{\rm a}$School of Mathematics and Statistics, Xi'an Jiaotong University, Xi'an 710049, P. R. China}\\
				{\footnotesize$^{\rm b}$School of Mathematics, Hunan University, Changsha 410082, P. R. China}
	}
	\date{}
	\maketitle
	\begin{abstract}
		Spectral characterizations of graphs is an important topic in spectral graph theory which has been studied extensively by researchers in recent years. The study of oriented graphs, however, has received less attention so far. In Qiu et al.~\cite{QWW} (Linear Algebra Appl. 622 (2021) 316-332), the authors gave an arithmetic criterion for an oriented graph to be determined by its \emph{generalized skew spectrum} (DGSS for short). More precisely, let $\Sigma$ be an $n$-vertex oriented graph with skew adjacency matrix $S$ and $W(\Sigma)=[e,Se,\ldots,S^{n-1}e]$ be the \emph{walk-matrix} of $\Sigma$, where $e$ is the all-one vector. A theorem of Qiu et al.~\cite{QWW} shows that a self-converse oriented graph $\Sigma$ is DGSS, provided that the Smith normal form of $W(\Sigma)$ is ${\rm diag}(1,\ldots,1,2,\ldots,2,2d)$, where $d$ is an odd and square-free integer and the number of $1$'s appeared in the diagonal is precisely $\lceil \frac{n}{2}\rceil$. In this paper, we show that the above square-freeness assumptions on $d$ can actually be removed, which significantly improves upon the above theorem. Our new ingredient is a key intermediate result, which is of independent interest: for a self-converse oriented graphs $\Sigma$ and an odd prime $p$, if the rank of $W(\Sigma)$ is $n-1$ over $\mathbb{F}_p$, then the kernel of $W(\Sigma)^{\rm T}$ over $\mathbb{F}_p$ is \emph{anisotropic}, i.e., $v^{\rm T}v\neq 0$ for any $0\ne v\in{{\rm ker}\,W(\Sigma)^{\rm T}}$ over $\mathbb{F}_p$.
	\end{abstract}

\section{Introduction}

Let $G$ be a simple graph with vertex set $V(G)=\{v_1,v_2,\ldots,v_n\}$ and edge set $E(G)$, where $n=|V(G)|$ is the order of the graph. The \emph{adjacency matrix} of the graph $G$ is the $n\times n$ matrix $A(G)=(a_{ij})$, where $a_{ij}=1$ if $v_{i}$ and $v_{j}$ are adjacent, and $a_{ij}=0$ otherwise. The \emph{adjacency spectrum} of $G$ consists of all the eigenvalues (including the multiplicities) of $A(G)$.

Two graphs are \emph{cospectral} if they share the same spectrum. For a real number $t\in\mathbb{R}$, two graphs $G$ and $H$ are \emph{$t$-cospectral} if $tJ-A(G)$ and $tJ-A(H)$ has the same spectrum, i.e., $\det(xI-(tJ-A(G)))= \det(xI-(tJ-A(H)))$, where $I$ is the identity matrix and $J$ is the all-one matrix. We say two graphs $G$ and $H$ are $\mathbb{R}$-\emph{cospectral} if they are $t$-cospectral for every $t\in \mathbb{R}$. A notable result of Johnson and Newman~\cite{JN} says that two graphs are $\mathbb{R}$-cospectral if and only if they are $t$-cospectral for any two distinct values of $t$.

 A graph $G$ is said to be \emph{determined by the spectrum} (DS for short), if any graph $H$ having the same spectrum as $G$ is isomorphic to $G$. A long standing unsolved question in spectral graph theory asks ``Which graphs are DS?''. The question originates from chemistry in 1956; see G\"{u}nthard and Primas~\cite{GP}. It is also closely related to other problems of central interest such as the graph isomorphism problem and the famous problem of Kac~\cite{KAC}: ``Can one hear the shape of a drum?". It is generally very hard and challenging to show a graph to be DS and up to now, very few families of graphs with special structures were shown to be DS. For more background and known results on this topic, we refer the reader to \cite{DH,DH1} and the references therein.

   Recently, Wang and Xu~\cite{W1,W2} and Wang~\cite{W3,W4} considered the above problem in the context of the generalized adjacency spectrum. We say two graphs $G$ and $H$ have the same \emph{generalized spectrum} if they are cospectral with cospectral complements, or equivalently, they are $\mathbb{R}$-cospectral by Johnson and Newman~\cite{JN}. A graph $G$ is \emph{determined by the generalized spectrum} (DGS for short), if any graph $H$ having the same generalized spectrum as $G$ is isomorphic to $G$.

Let $W(G)=[e,A(G)e,\ldots,A(G)^{n-1}e]$ be the \emph{walk-matrix} of graph $G$. The following theorem gives a simple arithmetic criterion for a graph to be DGS.
\begin{theorem}[Wang~\cite{W4}]\label{mainforgraph}
Let $G$ be graph of order $n$. If $2^{-\lfloor\frac{n}{2}\rfloor}\det W(G)$ (which is always an integer) is odd and square-free, then $G$ is DGS.
\end{theorem}

 The problem of generalized spectral characterizations naturally extends to oriented graphs.
Given a simple graph $G$ with an orientation $\sigma$. An \emph{oriented graph} $\Sigma=(G,\sigma)$ is a digraph obtained from $G$ by
assigning every edge a direction by $\sigma$, and $G$ is called the \emph{underlying graph} of $\Sigma$.

An oriented graph can be represented by its skew adjacency matrix, which was proposed by Tutte~\cite{Tu}, as a tool to deal with the problem of counting the number of perfect matchings in a graph. Let $S=S(\Sigma)=(S_{ij})$ be the \emph{skew adjacency matrix} of $\Sigma$ . Then
    \begin{equation*}
S_{i,j}= \begin{cases}
   1  &\mbox{if $(v_i,v_j)$ is an edge};\\
      -1 &\mbox{if $(v_j,v_i)$ is an edge};\\
   0  &\mbox{otherwise}.
   \end{cases}
\end{equation*}

\begin{example} Let $\Sigma$ be an oriented graph given as in Fig.~1.
\begin{figure}[h]
				\centering
				\includegraphics[width=.35\textwidth]{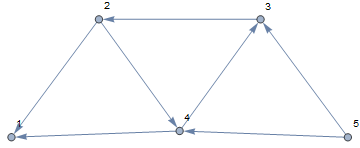}
\caption{An oriented graph $\Sigma$ with skew adjacency matrix $S$.}
								\label{img}
			\end{figure}
Then the skew adjacency matrix of $\Sigma$ is
 $$S=\left(
\begin{array}{ccccc}
 0 & -1 & 0 & -1 & 0 \\
 1 & 0 & -1 & 1 & 0 \\
 0 & 1 & 0 & -1 & -1 \\
 1 & -1 & 1 & 0 & -1 \\
 0 & 0 & 1 & 1 & 0 \\
\end{array}
\right).$$

\end{example}
The generalized spectral determination of graphs naturally extends to oriented graphs.
\begin{definition}
An oriented graph $\Sigma$ is said be determined by its generalized skew spectrum (DGSS for short), if whenever $\Delta$ is an oriented graph
such that $${\rm Spec}(S(\Sigma))={\rm Spec}(S(\Delta))~{\rm and}~{\rm Spec}(J-S(\Sigma))={\rm Spec}(J-S(\Delta)),$$
then $\Delta$ is isomorphic to $\Sigma$.

\end{definition}

The \emph{converse} of $\Sigma$, denoted by $\Sigma^{\rm T}$, is an oriented graph obtained by reversing every edge of $\Sigma$, and $\Sigma$ is called \emph{self-conversed} if it is isomorphic to its converse $\Sigma^{\rm T}$.
That is, $\Sigma$ is self-converse if there is a 1-1 mapping $\Phi: V(\Sigma)\rightarrow V(\Sigma)$ such that
$$(u,v)\in E(\Sigma)\iff(\Phi(v),\Phi(u))\in E(\Sigma).$$
The above $\Phi$ is also called an \emph{anti-automorphism} of $\Sigma$.
In other words, $\Sigma$ is self-conversed if and only if there exists a permutation matrix $P$ such that
\begin{equation}\label{Per}
P^{\rm T}SP=S^{\rm T} =-S.
\end{equation}

Let $\Phi$ be an anti-automorphism of $\Sigma$. Then $\Phi^2$ is an automorphism of $\Sigma$. Moreover, if $\Sigma$ is controllable (i.e., $\det W\neq0$), then $\Sigma$ has only one trivial automorphism - the identity mapping; see Lemma~\ref{dc} in Section 2.2. Thus, $\Phi^2=Id$, i.e., $\Phi$ is an involution or equivalently, the corresponding permutation matrix $P$ is symmetric and satisfies $P^2=I_n$. This property is crucial in our paper.
\begin{example} \textup{Let $\Sigma$ be an self-conversed oriented graph given as in Fig.~2.}
\begin{figure}[h]
				\centering
				\includegraphics[width=.30\textwidth]{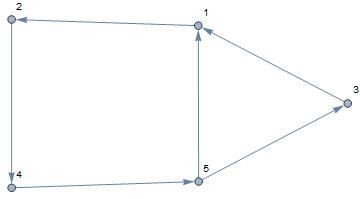}
				\caption{A self-converse oriented graph $\Sigma$} 
				\label{img}
			\end{figure}

$$S=\left(
\begin{array}{ccccc}
 0 & 1 & -1 & 0 & -1 \\
 -1 & 0 & 0 & 1 & 0 \\
 1 & 0 & 0 & 0 & -1 \\
 0 & -1 & 0 & 0 & 1 \\
 1 & 0 & 1 & -1 & 0 \\
\end{array}
\right),~P=\left(
\begin{array}{ccccc}
 0 & 0 & 0 & 0 & 1 \\
 0 & 0 & 0 & 1 & 0 \\
 0 & 0 & 1 & 0 & 0 \\
 0 & 1 & 0 & 0 & 0 \\
 1 & 0 & 0 & 0 & 0 \\
\end{array}
\right).$$

\end{example}

\begin{remark} \textup{We remark that the self-converseness assumption is necessary for $\Sigma$ to be DGSS. This is because $S(\Sigma^{\rm T})=S(\Sigma)^{\rm T}=-S(\Sigma)$, $\Sigma$ and $\Sigma^{\rm T}$ always have the same generalized skew spectrum, i.e., $\det(xI-tJ+S(\Sigma))=\det(xI-tJ+S(\Sigma)^{\rm T}))$ for every $t\in \mathbb{R}$. Thus, if an oriented graph $\Sigma$ is DGSS then it has to be self-converse.}
\end{remark}

Let $W({\Sigma}):=[e,Se,\ldots,S^{n-1}e]$ ($e$ is the all-one vector) be the \emph{walk-matrix} of $\Sigma$. Qiu et al.~\cite{QWW} proved the following theorem.

\begin{theorem}[Qiu et al.~\cite{QWW}]\label{QWW}
 Let ${\Sigma}$ be a self-converse oriented graph on $n$ vertices. Suppose that $2^{-\lfloor n/2\rfloor}\det W(\Sigma)$ is odd and square-free. Then $\Sigma$ is DGSS.
\end{theorem}

The condition in the above theorem can also be described in terms of the Smith normal form (SNF for short; see Section 2.1 for the definition) of $W(\Sigma)$. Actually, ${2^{-\lfloor n/2 \rfloor}}{\det W(\Sigma)}$ is odd and square-free if and only if the SNF of $W(\Sigma)$ has the following form:
\begin{equation}\label{PP1}
{\rm diag}(\underbrace{1,1,\ldots,1}_{\lceil\frac{n}2\rceil},\underbrace{2,2,\ldots,2,2d}_{\lfloor \frac{n}2 \rfloor}),
\end{equation}
where $d$ is an odd and square-free integer.
That is, a self-conversed oriented graph $\Sigma$ is DGSS provided Eq.~\eqref{PP1} holds.

The main result of the paper is the following theorem, which shows that the square-freeness assumption of $d$ in Eq.~\eqref{PP1} can actually be removed, and thus significantly improves upon Theorem~\ref{QWW}.

\begin{theorem}\label{main}Let $\Sigma$ be a self-converse oriented graph with walk-matrix $W(\Sigma)$.
Suppose that the SNF of $W(\Sigma)$ is as follows:
$${\rm diag}(\underbrace{1,1,\ldots,1}_{\lceil n/2 \rceil},\underbrace{2,2,\ldots,2,2d}_{\lfloor n/2\rfloor}),$$
where $d$ is any odd integer. Then $G$ is DGSS.
\end{theorem}

We remark that the proof of Theorem~\ref{QWW} is based on the methods developed in \cite{W1,W3,W4}. While the proof of Theorem~\ref{main} is based on a new discovery which greatly simplifies and strengthens that of Theorem~\ref{QWW}. To sate the key idea, we need a definition below.

\begin{definition}
\textup{Let $p$ be a prime. A subspace $U\subset \mathbb{F}_p^n$ is \emph{totally isotropic} if $u^{\rm T}w=0$ for any $u,w\in{U}$,
and $U\subset \mathbb{F}_p^n$ is \emph{anisotropic} if $u^{\rm T}u\neq 0$ for any $0\neq u\in{U}$.}
\end{definition}

Central to the proof of Theorem~\ref{main} is the following theorem, which is of independent interest.

\begin{theorem}\label{tournament}Let $\Sigma$ be a self-conversed oriented graph on $n$ vertices such that $\det W(\Sigma)\ne 0$.
 Let $p$ be an odd prime. Suppose that ${\rm rank}_p\,W(\Sigma)=n-1$, where ${\rm rank}_p\,W(\Sigma)$ denotes the rank of $W(\Sigma)$ over $\mathbb{F}_p$. Then $\ker W(\Sigma)^{\rm T}\subset \mathbb{F}_p^n$ is anisotropic, i.e.,
if $W(\Sigma)^{\rm T}v=0$ with $v\neq 0$, then $v^{\rm T}v\neq 0$ over $\mathbb{F}_p$.
\end{theorem}

We shall postpone the proofs of Theorem~\ref{main} and \ref{tournament} to Section 3. We remark that results similar to Theorem~\ref{main} and \ref{tournament} are not available for simple graphs. Thus, the problem of spectral characterization of oriented graphs has some particular characteristics, and deserves an exclusive treatment in its own right.

The rest of the paper is organized as follows. In Section 2, we give some preliminary results that will be needed later in the paper. In Section 3, we present the proof of Theorems~\ref{main} and \ref{tournament}. Some examples are given in Section 4 to illustrate Theorem~\ref{main}.

\noindent
\textbf{Notations}: In the rest of the paper, we shall write $S=S(\Sigma)$ and $W=W(\Sigma)$ if there is no confusion arises. Throughout, let $p$ be an odd prime and $\mathbb{F}_p$ is the finite filed with $p$ elements. For an integral matrix $M$, we shall use ${\rm rank}_p M$ and ${\rm rank}\,M$ to denote the rank of $M$ over $\mathbb{F}_p$ and $\mathbb{Q}$, respectively.

\section{Preliminaries}

For the convenience of the reader, we shall give some preliminary results that will be needed later in the paper.

\subsection{The main strategy}

First, we shall describe our main strategy for showing an oriented graph to be DGSS, which roughly follows the general framework in~\cite{W1,W3,W4}. The following lemma gives a simple characterization of two oriented graphs having the same generalized skew spectrum. The proof is omitted here since the previous proof in~\cite{W1} also applies to the generalized skew spectrum with some slight modifications.

\begin{lemma}[c.f. Wang and Xu~\cite{W1}] \label{rational} Let $\Sigma$ be an oriented graph with $\det W(\Sigma)\neq 0$. There exists an oriented graph $\Delta$ such that $\Sigma$ and $\Delta$ have the same generalized skew spectrum if and only if there exists a unique rational orthogonal matrix $Q$ such that
\begin{equation}\label{EQ1}
Q^{\rm T}S(\Sigma)Q=S(\Delta)~{\rm and}~Qe=e.
\end{equation}
\end{lemma}

Define $$\mathscr{Q}(\Sigma)=\{Q\in{O_n(\mathbb{Q})}:\,Q^{\rm T}S(\Sigma)Q=S(\Delta)~{\rm for~some~graph}~\Delta~{\rm and}~Qe=e\},$$
 where $O_n(\mathbb{Q})$ denotes the set of all orthogonal matrices with rational entries.

\begin{lemma}[Wang and Xu~\cite{W1}; Qiu et al.~\cite{QWW}]\label{good} Let $\Sigma$ be an oriented graph with $\det W(\Sigma)\neq 0$. Then $\Sigma$ is determined by its generalized skew spectrum if and only if $\mathscr{Q}(\Sigma)$ contains only permutation matrices.
\end{lemma}

By the above theorem, next, in order to show a given oriented graph $\Sigma$ is DGSS, we have to determine whether $\mathscr{Q}$$(\Sigma)$ contains only permutation matrices. In order to do so, we give the following definition.

\begin{definition}[Wang and Xu~\cite{W1}] Let $Q$ be an orthogonal matrix with rational entries. The level of $Q$, denoted by $\ell(Q)$ or simply $\ell$, is the smallest positive integer $k$ such that $kQ$ is an integral matrix.
\end{definition}

Clearly, a rational orthogonal matrix $Q$ with $Qe=e$ is a permutation matrix if and only if $\ell(Q)=1$. Thus, for a given oriented graph ${\Sigma}$, our main strategy in proving $\mathscr{Q}(\Sigma)$ contains only permutation matrices is to show that every $Q\in {\mathscr{Q}(\Sigma)}$ has level $\ell=1$.

When dealing with integral and rational matrices, the Smith normal form (SNF for short) is a useful tool. An integral matrix $U$ is called \emph{unimodular}, if $\det U=\pm1$.
The following theorem is well-known.

 \begin{theorem}[e.g.,~\cite{SC}] Let $M$ be an $n$ by $n$ integral matrix with full rank. Then there exist two unimodular matrices $U_1$ and $U_2$ such that $M=U_1NU_2$, where
$N={\rm diag}(d_1,d_2,\ldots,d_n)$ with $d_i\mid d_{i+1}$ for $i=1,2,\ldots,n-1$.
\end{theorem}
The above diagonal matrix $N$ is known as the \textit{Smith Normal Form} (SNF) of $M$, and
 $d_i$ is the $i$-th \emph{invariant factor} of $M$.

 The following lemma shows the level of a $Q\in {\mathscr{Q}(\Sigma)}$ is closely related to the $n$-th invariant factor of $W$.

\begin{lemma}[c.f.~Wang and Xu~\cite{W1}]\label{level} Let $\Sigma$ be an oriented graph. Let $Q\in {\mathscr{Q}(\Sigma)}$ with level $\ell$. Then $\ell\mid d_n$, where
$d_n$ is the $n$-th invariant factor of $W(\Sigma)$.
\end{lemma}

\subsection{Some lemmas}

In this subsection, we collect some lemmas needed in the proof of Theorems~\ref{main} and \ref{tournament}.

\begin{lemma}[Qiu et al.~\cite{QWW}]\label{even}Suppose the conditions of the Theorem~\ref{main} hold. Let $Q\in{\mathscr{Q}(\Sigma)}$ with level $\ell$. Then $\ell$ is odd.
\end{lemma}

\begin{lemma}[Qiu et al.~\cite{QWW}] Suppose that ${\rm rank}_p\,{W}=r$. Then the first $r$ columns of $W$ consists of a basis of the column space of $W$.

\end{lemma}

\begin{lemma}[Qiu et al.~\cite{QWW}]\label{ws}Let $\Sigma$ be a self-conversed oriented graph. Suppose that ${\rm rank}_p\, W=n-1$. If $W^{\rm T} v=0$ over $\mathbb{F}_p$, then $Sv=0$ over $\mathbb{F}_p$.
\end{lemma}
\begin{proof} We give a proof for completeness. It follows from $W^{\rm T}v\equiv 0~({\rm mod}~p)$ that $e^{\rm T}(S^{\rm T})^kv\equiv 0~({\rm mod}~p)$ for $k=0,1,\ldots,n-1$.
By the Cayley-Hamilton Theorem we have $S^ne=-(a_1S^{n-1}e+\cdots+a_nI)$. Thus, $e^{\rm T}(S^{\rm T})^nv\equiv 0~({\rm mod}~p)$. So we
have $W^{\rm T}(Sv)\equiv 0~({\rm mod}~p)$. Then the assumption that ${\rm rank}_p\, W=n-1$ implies that
\begin{equation}\label{E111}
Sv\equiv \lambda_0 v~({\rm mod}~p),
\end{equation}
for some $\lambda_0\in {\mathbb{Z}}.$

\noindent
\textbf{Claim:} $ \lambda_0 \equiv 0~({\rm mod}~p)$.

Actually, note that $P^{\rm T}SP=S^{\rm T}=-S$. It follows from Eq.~\eqref{E111} that $S(Pv)\equiv-\lambda_0(Pv)~({\rm mod}~p)$. Thus, we have $e^{\rm T}S^i(Pv)\equiv (-1)^i\lambda_0^ie^{\rm T}(Pv)~({\rm mod}~p)$ for $0\leq i\leq n-1$. Note that $e^{\rm T}(Pv)\equiv e^{\rm T}v\equiv 0~({\rm mod}~p)$. Thus we have $W^{\rm T}(Pv)\equiv 0~({\rm mod}~p)$. Further note that ${\rm rank}_p\,{W}=n-1$. So $Pv$ and $v$ are linearly dependent. Combing $Sv\equiv-\lambda_0v~({\rm mod}~p)$ and $Sv\equiv\lambda_0v~({\rm mod}~p)$ gives that $ \lambda_0 \equiv 0~({\rm mod}~p)$, as desired.
\end{proof}

\begin{lemma}\label{pfy} Let $S^{\rm T}=P^{\rm T}SP$. Then $PS^ke=(-1)^kS^ke$ for any $k\geq 0$.
\end{lemma} 
\begin{proof}
This follows directly by noticing $(P^{\rm T}SP)^ke=P^{\rm T}S^kPe=(-S)^ke$ and $Pe=e$.
\end{proof}

\begin{lemma}\label{dc}Let $Q$ be a regular orthogonal matrix such that $Q^{\rm T}SQ=S^{\rm T}$. Suppose that $W$ is non-singular. Then $Q$ is symmetric.
In particular, let $P$ be the permutation matrix such that $P^{\rm T}SP=S^{\rm T}$, then $P$ is symmetric and $P^2=I_n$.
\end{lemma}
\begin{proof}
It follows from $Q^{\rm T}SQ=S^{\rm T}$ that $S=QS^{\rm T}Q^{\rm T}$, i.e., $-S=S^{\rm T}=QSQ^{\rm T}$, or equivalently,
$(Q^{\rm T})^{\rm T}S(Q^{\rm T})=S^{\rm T}$. Thus, we have $Q^{\rm T}W(\Sigma)=W(\Sigma^{\rm T})$ and $QW(\Sigma)=W(\Sigma^{\rm T})$. It follows that
$Q^{\rm T}=Q=W(\Sigma^{\rm T})W(\Sigma)^{-1}$.
\end{proof}


\begin{lemma} Let $P$ be the permutation matrix in Eq.~\eqref{Per}. Then for even $n$, we have ${\rm rank}(I-P)={\rm rank}(I+P)=n/2$ and for odd $n$, we have ${\rm rank}(I-P)=(n-1)/2$ and ${\rm rank}(I+P)=(n+1)/2$. Similar results hold over $\mathbb{F}_p$.
\end{lemma}

\begin{proof} We only prove the even case. Over $\mathbb{Q}$, for odd $i$, we have $(I+P)S^i e=0$, i.e., $S^ie\in \ker (I+P)$. Since $\{Se,S^3e,\ldots,S^{n-1}e\}$ are linearly independent $n/2$ vectors, we have ${\rm rank}(I+P)\leq n/2$. Similarly, for even $i$, ${\rm rank}(I-P)\leq n/2$. Note
$$n\leq{\rm rank}_p(I+P)+{\rm rank}_p(I-P)\leq{\rm rank} (I+P)+{\rm rank} (I-P)\leq\frac n2+\frac n2=n,$$
so all inequalities become equalities and the lemma follows.
\end{proof}


Let $V_\lambda$ denote the eigenspace of $P$ corresponding to the eigenvalue $\lambda$ over $\mathbb{F}_p$. Then we have

\begin{corollary} Let $P$ be the permutation matrix such that $P^{\rm T}SP=S^{\rm T}$. Then for even $n$, $\dim V_1=\dim V_{-1}=\frac{n}2;$ for odd $n$,
$\dim V_1=\frac{n+1}2$ and $\dim V_{-1}=\frac{n-1}2.$

\end{corollary}


\section{Proofs of Theorem~\ref{main} and Theorem~\ref{tournament}}

In this section, we present the proofs of Theorem~\ref{main} and Theorem~\ref{tournament}.

Notice that under the assumptions of Theorem~\ref{main}, we have $r={\rm rank}_p\,{W}=n-1$ and the first $n-1$ columns of $W$ are linearly independent over $\mathbb{F}_p$. The following notion turns out to be useful in proving our main theorem.

\begin{definition}
Define $$V_0:={\rm Col}_{\mathbb{F}_p} \,(W)={\rm span}\langle e,Se,\ldots,S^{n-2}e \rangle\subset \mathbb{F}_p^n.$$
\end{definition}


Now we present the proof of Theorem~\ref{tournament}.

\begin{proof}[Proof of Theorem~\ref{tournament}]

 We only prove the case that $n$ is even, the case $n$ is odd can be proved similarly.
 Since $\rank_p W=n-1$, we have $\dim V_0^\perp =1$ and $V_0^\perp={\rm span}\langle v\rangle$.
 By Lemma~\ref{pfy}, $PS^{i}e=(-1)^iS^{i}e$. It follows that $e,S^2e,\dots,S^{n-2}e$ are the $\frac{n}2$ eigenvalues of $P$ corresponding to $1$, and $Se,S^3e,\dots,S^{n-3}e$ are the $\frac{n}2-1$ eigenvalues of $P$ corresponding to $-1$. Note the dimension of the eigenspace of $P$ corresponding to $-1$ is equal to $\frac{n}2$. Thus, there exists an $\gamma \in V\setminus V_0$ such that $P\gamma=-\gamma$.

 For contradiction, suppose that $v^{\rm T}v=0$. Then $v^{\rm T}[W,v]=0$. Note that $\rank_p\,W=n-1$. It follows that $v\in {V_0}$.
 Let $v=\sum_{i=0}^{n-2}c_{i}S^ie$. Then we have

 \noindent
\textbf{ Claim:} $c_{n-2}\neq 0$.

Otherwise, suppose $v=\sum_{i=0}^{n-3} c_i S^i e$. Since $W^{\rm T} v=0$, by Lemma~\ref{ws}, we have $Sv=0$, which implies that $\sum_{i=1}^{n-2} c_{i-1} S^i e=0$. Since $Se,S^2 e,\ldots, S^{n-2}e$ are linearly independent, we have $c_i=0$ for all $i$, contradicting $v\neq 0$.

Note that $PS^{i}e=(-1)^iS^{i}e$, we have $Pv\in V_0^\perp$, so $Pv=kv$ for some $k\in {\mathbb{F}_p}$. Comparing the coefficients of $S^{n-2}e$ on both sides gives that $k=1$, i.e., $Pv=v$.
     So we have
$$-v^{\rm T}\gamma=v^{\rm T}P\gamma=v^{\rm T}\gamma.$$
 It follows that $v^{\rm T}\gamma=0$, and hence $v$ is orthogonal to every vector in $\mathbb{F}_p^n$ since $\mathbb{F}_p^n=V_0\bigoplus {\rm span}\langle \gamma\rangle$; a contradiction! Therefore we have $v^{\rm T}v\ne0$ and the theorem follows.

\end{proof}

Finally, we are ready to present the proof of Theorem~\ref{main}.

\begin{proof}[Proof of Theorem~\ref{main}] Let $\Delta$ be any oriented graph that is generalized cospectral with $\Sigma$. Then there exists a $Q\in{\mathscr{Q}(\Sigma)}$ with level $\ell$ such that $Q^{\rm T}S(\Sigma)Q=S(\Delta)$. It follows $Q^{\rm T}W(\Sigma)=W(\Delta)$.
For contradiction, suppose that $\ell\neq 1$. Let $p$ be any prime divisor of $\ell$. Then it follows from Lemma~\ref{level} that $\ell\mid 2d$. According to Lemma~\ref{even}, $\ell$ is odd. Thus, $p$ must be odd. Moreover, we have $p\mid d$ and ${\rm rank}_p\,{W(\Sigma)}=n-1$.
 Let $U:={\rm Col}_{\mathbb{F}_p}({\ell Q})$.

 \noindent
\textbf{ Claim:}  $U\subset \ker W(\Sigma)^{\rm T}$ is totally isotropic.

Let $u$ be any column of $\ell Q$. Then it follows from  $W(\Sigma)^{\rm T}(\ell Q)=\ell W(\Delta)^{\rm T}\equiv 0~({\rm mod}~p)$ that $W(\Sigma)^{\rm T}u\equiv 0~({\rm mod}~p)$, and hence $U\subset \ker W(\Sigma)^{\rm T}$. Let $u,w$ be any two columns of $\ell Q$, then either $u^{\rm T}w=\ell^2$ or $0$, which implies that $U$ is totally isotropic.

     Nevertheless, this contradicts Theorem~\ref{tournament}. Thus, $\ell=1$ and $Q$ is a permutation matrix. Hence $\Delta$ is isomorphic to $\Sigma$. This completes the proof.
\end{proof}

\section{Some examples}
In this section, we shall give some examples for an illustration of Theorem~\ref{main}.

\begin{example} \textup{Let $\Sigma$ be a self-conversed oriented graph given as in Fig~3.}
\begin{figure}[h]
				\centering
				\includegraphics[width=.35\textwidth]{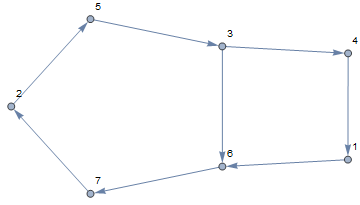}
				\caption{A self-conversed DGSS oriented graph $\Sigma$} 
				\label{img}
			\end{figure}
$$S=\left(
\begin{array}{ccccccc}
 0 & 0 & 0 & -1 & 0 & 1 & 0 \\
 0 & 0 & 0 & 0 & 1 & 0 & -1 \\
 0 & 0 & 0 & 1 & -1 & 1 & 0 \\
 1 & 0 & -1 & 0 & 0 & 0 & 0 \\
 0 & -1 & 1 & 0 & 0 & 0 & 0 \\
 -1 & 0 & -1 & 0 & 0 & 0 & 1 \\
 0 & 1 & 0 & 0 & 0 & -1 & 0 \\
\end{array}
\right),~P=\left(
\begin{array}{ccccccc}
 0 & 0 & 0 & 1 & 0 & 0 & 0 \\
 0 & 1 & 0 & 0 & 0 & 0 & 0 \\
 0 & 0 & 0 & 0 & 0 & 1 & 0 \\
 1 & 0 & 0 & 0 & 0 & 0 & 0 \\
 0 & 0 & 0 & 0 & 0 & 0 & 1 \\
 0 & 0 & 1 & 0 & 0 & 0 & 0 \\
 0 & 0 & 0 & 0 & 1 & 0 & 0 \\
\end{array}
\right).$$

\textup{The SNF of $W$ is ${\rm diag}(1,1,1,1,2,2,2\times3^2)$. Thus, according to Theorem~\ref{main}, $\Sigma$ is DGSS.}

\end{example}

\begin{example} {\textup {Let $\Sigma$ be a self-conversed oriented graph given as in Fig~4.}}

\begin{figure}[h]
				\centering
				\includegraphics[width=.35\textwidth]{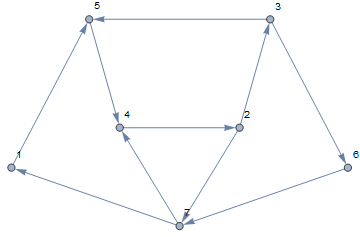}
				\caption{A self-conversed DGSS oriented graph $\Sigma$} 
				\label{img}
			\end{figure}
$$S=\left(
\begin{array}{ccccccc}
 0 & 0 & 0 & 0 & 1 & 0 & -1 \\
 0 & 0 & 1 & -1 & 0 & 0 & 1 \\
 0 & -1 & 0 & 0 & 1 & 1 & 0 \\
 0 & 1 & 0 & 0 & -1 & 0 & -1 \\
 -1 & 0 & -1 & 1 & 0 & 0 & 0 \\
 0 & 0 & -1 & 0 & 0 & 0 & 1 \\
 1 & -1 & 0 & 1 & 0 & -1 & 0 \\
\end{array}
\right),~
P=\left(
\begin{array}{ccccccc}
 0 & 0 & 0 & 0 & 0 & 1 & 0 \\
 0 & 0 & 0 & 1 & 0 & 0 & 0 \\
 0 & 0 & 0 & 0 & 1 & 0 & 0 \\
 0 & 1 & 0 & 0 & 0 & 0 & 0 \\
 0 & 0 & 1 & 0 & 0 & 0 & 0 \\
 1 & 0 & 0 & 0 & 0 & 0 & 0 \\
 0 & 0 & 0 & 0 & 0 & 0 & 1 \\
\end{array}
\right).$$

{\textup {The SNF of $W$ is ${\rm diag}(1,1,1,1,2,2,2\times3^2)$. Thus, according to Theorem~\ref{main}, $\Sigma$ is DGSS.}}

\end{example}

\section*{Acknowledgments}
The research of the second author is supported by National Key Research and Development Program of China 2023YFA1010203 and National Natural Science Foundation of China (Grant No.\,12371357), and the third author is supported by Fundamental Research Funds for the Central Universities (Grant No.\,531118010622), National
Natural Science Foundation of China (Grant No.\,1240011979) and Hunan Provincial Natural Science Foundation of China (Grant No.\,2024JJ6120).


\begin{thebibliography}{22}




\bibitem{DH}
E.R. van Dam, W.H. Haemers, Which graphs are determined by their spectrum?
\emph{Linear Algebra Appl.} 373 (2003) 241--272.

\bibitem{DH1}
 E.R. van Dam, W.H. Haemers, Developments on spectral characterizations
of graphs, \emph{Discrete Math.} 309 (2009) 576--586.

\bibitem{GP} Hs.H. G\"{u}nthard, H. Primas, Zusammenhang von Graphentheorie und MO-Theorie von Molekeln
mit Systemen konjugierter Bindungen, \emph{Helv. Chim. Acta,} 39 (1956) 1645--1653.

\bibitem{HK} K. Hoffman, R. Kunze, \emph{Linear Algebra}, (2nd edition), PHI, 2009.




\bibitem{JN} C.R. Johnson, M.~Newman, A note on cospectral graphs, \emph{J. Combin. Theory, Ser. B,} 28(1980) 96--103.

\bibitem{KAC}  M. Kac, Can one hear the shape of a drum? \emph{Amer. Math. Monthly},  73(4) (1966) 1-23.

\bibitem{QWW} L.H. Qiu, W. Wang, W. Wang, On oriented graphs determined by their generalized skew spectra, \emph{Linear Algebra Appl}. 622 (2021) 316-332.


\bibitem{QWWZ} L.H. Qiu, W. Wang, W. Wang, H. Zhang, Smith normal form and the generalized spectral characterization of graphs,
\emph{Discrete Math}. 346 (2023) 113177.

\bibitem{SC}
A. Schrijver, \emph{The Theory of Linear and Integer Programming}, John Wiley $\&$ Sons, 1998.

\bibitem{Tu} W.T. Tutte, The factorization of linear graphs, \emph{J. London Math. Soc.,} 22 (1947) 107-111.



\bibitem{W1}
 W. Wang, C. X. Xu, A sufficient condition for a family of graphs being determined
by their generalized spectra, \emph{European J. Combin.}, 27 (2006) 826--840.

\bibitem{W2}
 W. Wang, C.X. Xu, An excluding algorithm for testing whether a family of
graphs are determined by their generalized spectra, \emph{Linear Algebra and its
Appl.}, 418 (2006) 62--74.


\bibitem{W3} W. Wang, Generalized spectral characterization of graphs revisited, \emph{Electron. J.
Comb.}, 20 (4) (2013), $\sharp$ P4.


\bibitem{W4} W. Wang, A simple arithmetic criterion for graphs being
determined by their generalized spectra, \emph{J. Combin. Theory, Ser. B}, 122 (2017) 438--451.





\end{thebibliography}
\end{document}